\documentclass{amsart}

\usepackage{amssymb}
\usepackage{amscd}
\usepackage{latexsym}
\usepackage{amsfonts}
\usepackage{amsmath}

\newtheorem{proposition}{Proposition}
\newtheorem{lemma}{Lemma}
\newtheorem{theorem}{Theorem}
\newtheorem{corollary}{Corollary}
\newtheorem{definition}{Definition}

\newcommand{\eproof}{\begin{flushright} $\square$ \end{flushright}}

\newcommand{\tr}{\mathop{\fam0 Tr}\nolimits}

\newcommand{\rank}{\mathop{\fam0 Rank}\nolimits}

\newcommand{\Hom}{\mathop{\fam0 Hom}\nolimits}
\newcommand{\End}{\mathop{\fam0 End}\nolimits}

\newcommand{\ad}{\mathop{\fam0 ad}\nolimits}

\newcommand{\id}{\mathop{\fam0 Id}\nolimits}

\newcommand{\Aut}{\mathop{\fam0 Aut}\nolimits}

\newcommand{\Pic}{\mathop{\fam0 Pic}\nolimits}

\newcommand{\bC}{{\mathbb C}}

\newcommand{\C}{C}

\newcommand{\Z}{{\mathbb Z}}
\newcommand{\bZ}{\Z{}}

\newcommand{\ra}{\mathop{\fam0 \rightarrow}\nolimits}

\newcommand{\M}{{\mathcal M}}

\newcommand{\cH}{ {\mathcal H}}
\newcommand{\T}{ {\mathcal T}}
\newcommand{\fg}{ {\mathfrak g}}
\newcommand{\V}{ {\mathcal V}}
\renewcommand{\L}{{\mathcal L}}

\renewcommand{\P}{ {\mathbb P}}

\newcommand{\s}{\sigma}
\newcommand{\Nabla}{{\mathbf {\hat \nabla}}}
\newcommand{\Nablae}{{\mathbf {\hat \nabla}}^e}
\newcommand{\Nablaet}{{\mathbf {\hat \nabla}}^{e,t}}

\newcommand{\Teim}{Teichm{\"u}ller }
\newcommand{\Ric}{\mathop{\fam0 Ric}\nolimits}

\begin{document}

\title[Asymptotic faithfulness]{Asymptotic faithfulness of the quantum $SU(n)$
representations of the mapping class groups}

\author{J\o rgen Ellegaard Andersen}
\address{Department of Mathematics\\
        University of Aarhus\\
        DK-8000, Denmark}
\email{andersen@imf.au.dk}
 \thanks{This research was conducted
 for the Clay Mathematics Institute at University of California,
 Berkeley.}

\begin{abstract}
We prove that the sequence of projective quantum $SU(n)$ representations of the
mapping class group of a closed oriented surface, obtained from the projective flat $SU(n)$-Verlinde
bundles over \Teim space, is asymptotically
faithful, that is the intersection over all levels of the kernels of
these representations is trivial, whenever the genus is at least $3$. For the genus $2$
case, this intersection is exactly the order two subgroup, generated
by the hyper-elliptic involution, in the case of even degree and $n=2$. Otherwise the
intersection is also trivial in the genus $2$ case.
\end{abstract}

\maketitle

\section{Introduction}

In this paper we shall study the finite dimensional
{\em quantum $SU(n)$
representations} of the mapping class group of a genus $g$ surface.
These form the only rigorously constructed part of the gauge-theoretic approach to
Topological Quantum Field Theories
in dimension $3$, which Witten proposed in his seminal paper \cite{W1}.
We discovered the {\it asymptotic faithfulness}
property for these representations by studying this approach,
which we will now briefly describe, leaving further details to
section \ref{sec2}, \ref{sec3} and the references given there.

Let $\Sigma$ be a closed oriented surface of genus $g\geq 2$ and $p$ a point on
$\Sigma$. Fix  $d\in \Z/n\Z \cong Z_{SU(n)}$ in the center of $SU(n)$. Let $M$ be the moduli
space of flat $SU(n)$-connections on $\Sigma - p$ with holonomy
$d$ around $p$.

By applying geometric quantization to the moduli space $M$ one
gets a certain finite rank vector bundle over \Teim space $\T$, which we
will call the {\em Verlinde} bundle $\V_k$ at level $k$, where $k$
is any positive integer. The fiber of this bundle over a point
$\sigma\in \T$ is $\V_{k,\sigma} = H^0(M_{\sigma},\L_\s^k)$, where
$M_{\s}$ is $M$ equipped with a complex structure induced from
$\sigma$ and $\L_\s$ is an ample generator of the Picard group of
$M_{\s}$.

The main result pertaining to this bundle $\V_k$ is that its
projectivization $\P(\V_k)$ supports a natural flat connection.
This is a result proved independently by Axelrod, Della Pietra and
Witten \cite{ADW} and by Hitchin \cite{H}. Now, since there is an
action of the mapping class group $\Gamma$ of $\Sigma$ on $\V_k$
covering its action on $\T$, which preserves the flat connection
in $\P(\V_k)$, we get for each $k$ a finite dimensional projective
representation, say $\rho^{n,d}_k$, of $\Gamma$, namely on the
covariant constant sections of $\P(\V_k)$ over $\T$. This sequence
of projective representations $\rho^{n,d}_k$, $k\in {\mathbb N}_+$ is the {\em quantum $SU(n)$
representations} of the mapping class group $\Gamma$.

For each given $(n,d,k)$, we cannot expect $\rho^{n,d}_k$ to be
faithful. However, V. Turaev conjectured a decade ago (see e.g. \cite{T})
that there should be no non-trivial element of the mapping class group
representing trivially under $\rho^{n,d}_k$ for all $k$, keeping
$(n,d)$ fixed. We call this property {\em asymptotic faithfulness} of the
quantum $SU(n)$ representations $\rho^{n,d}_k$. In this paper we
prove Turaev's conjecture:

\begin{theorem}\label{Main} Assume that $n$ and $d$ are coprime or
  that $(n,d)=(2,0)$ when $g=2$.
Then we have that
\begin{equation*}
\bigcap_{k=1}^\infty \ker(\rho^{n,d}_k) =
\begin{cases}
\{1, H\} & g=2 \mbox{, }n=2 \mbox{ and } d=0 \\
\{1\}& \mbox{otherwise}.
\end{cases}
\end{equation*}
where $H$ is the
hyperelliptic involution.
\end{theorem}

This Theorem states that for any element $\phi$ of the mapping
class group $\Gamma$, which is not the identity element (and not the hyperelliptic involution
in genus 2), there is
an integer $k$ such that $\rho^{n,d}_k(\phi)$ is not a multiple of
the identity. We will suppress the superscript on the quantum
representations and simply write $\rho_k = \rho^{n,d}_k$
throughout the rest of the paper.

Our key idea in the proof of Theorem \ref{Main} is the use of the endomorphism
bundle $\End(\V_k)$ and the construction of sections of this
bundle via {\em Toeplitz operators} associated to smooth functions on
the moduli space $M$. By showing that these sections are {\em
asymptotically} flat sections of $\End(\V_k)$ (see Theorem
\ref{Asympflat} for the precise statement), we prove that any
element in the above intersection of kernels acts trivially on the
smooth functions on $M$, hence acts by the identity on $M$ (see the
proof of Theorem \ref{Main1}). Theorem \ref{Main} now follows
directly from knowing which elements of the mapping class group
act trivially on the moduli space $M$.

The assumptions on the pair $(n,d)$ in Theorem \ref{Main} are exactly picking out the
cases where the moduli space $M$ is smooth.
This means we can apply the work of Bordemann, Meinrenken
and Schlichenmaier on Toeplitz operators on smooth K{\"a}hler
manifolds, in particular their formula for the asymptotics
in $k$ of the operator norm of Toeplitz operators and the
asymptotic expansion of the product of two Toeplitz operators. -
It is with the use of these results we establish that the
Toeplitz operator sections are asymptotically flat with respect to
Hitchin's connection.

In the remaining cases, where $M$ is singular, we also have a proof of asymptotic
faithfulness, where we use the desingularization of the moduli
space, but this argument is technically quite a bit more involved.
However, together with Michael Christ
we have in \cite{AC} extended some of the results
of Bordemann, Meinrenken and Schlichenmaier and Karabegov and
Schlichenmaier to the case of singular varieties.
In \cite{A3} the argument of the present paper
will be repeated in the non-coprime case, making use of the results
of \cite{AC} to show that Theorem \ref{Main}
holds in general without the coprime assumption.

The abelian case, i.e. the case where $SU(n)$ is replaced by
$U(1)$, was considered in \cite{A2}, before we considered the case
discussed in this paper. In this case, with the use of
Theta-functions, explicit expressions for the Toeplitz operators
associated to holonomy functions can be obtained. From these
expressions it follows that the Toeplitz operators are not
covariant constant even in this much simpler case
(although the relevant connection is the one induced from the
$L_2$-projection as shown by Ramadas in \cite{R1}). However, they
are asymptotically covariant constant, in fact we find explicit
perturbations to all orders in $k$, which in this case, we argue,
can be summed and used to create actual covariant constant sections
of the endomorphism bundle. The result as far as the mapping class
group goes, is that the intersection of the kernels over all $k$
in that case is the Torelli group.

Returning to the non-abelian case at hand, we know by the
work of Laszlo \cite{La1}, that $\P(\V_k)$ with
its flat connection is isomorphic to the projectivization of the
bundle of conformal blocks for $sl(n,\bC)$ with its flat connection
over $\T$ as constructed by Tsuchiya, Ueno and Yamada \cite{TUY}.
This means that the quantum $SU(n)$ representations $\rho_k$ is
the same sequence of representations as the one arising from the
space of conformal blocks for $sl(n,\bC)$.

By the work of Reshetikhin-Turaev, Topological Quantum Field
Theories have been constructed in dimension $3$ from the quantum group $U_q sl(n,\bC)$
(see \cite{RT1}, \cite{RT2} and \cite{T})
or alternatively from the Kauffman bracket and the Homfly-polynomial
by Blanchet, Habegger, Masbaum and Vogel (see \cite{BHMV1}, \cite{BHMV2} and
\cite{B1}).

In ongoing work of Ueno joint with this author (see \cite{AU1},
\cite{AU2} and \cite{AU3}), we are in the process of establishing
that the TUY construction of the bundle of conformal blocks over
\Teim space for $sl(n,\bC)$ gives a modular functor, which in turn
gives a TQFT, which is isomorphic to the $U_q
sl(n,\bC)$-Reshetikhin-Turaev TQFT. A corollary of this will be that the
quantum $SU(n)$ representations are isomorphic to the ones that
are part of the $U_q sl(n,\bC)$-Reshetikhin-Turaev TQFT.
 Since it is well known that
the Reshetikhin-Turaev TQFT is unitary one will get unitarity of the quantum
$SU(n)$ representations from this. We note that unitarity is not clear from
the
geometric construction of the quantum $SU(n)$ representions.
If the quantum $SU(n)$ representations $\rho_k$ are unitary, then
we have for all $\phi\in \Gamma$ that
\begin{equation}
|\tr(\rho_k(\phi))| \leq \dim \rho_k.\label{ineqaul1}
\end{equation}

From the above Theorem \ref{Main} follows directly, assuming
unitarity, that one has the following
\begin{corollary}\label{cor1}
Assume that $n$ and $d$ are coprime or
  that $(n,d)=(2,0)$ when $g=2$.
Then equality holds in (\ref{ineqaul1}) for all $k$, if and only if
\begin{equation*}
\phi \in
\begin{cases}
\{1, H\} & g=2 \mbox{, }n=2 \mbox{ and } d = 0\\
\{1\}& \mbox{otherwise}.
\end{cases}
\end{equation*}
\end{corollary}

Furthermore, one will get that
the norm of the Reshetikhin-Turaev quantum invariant for all $k$ and $n=2$
($n=3$ in the genus $2$ case)
can separate the mapping
torus of the identity from the rest of the mapping tori as a purely TQFT
consequence of Corollary \ref{cor1}.

In this paper we have initiated the program of using of the theory of Toeplitz
operators on the moduli spaces in the study of TQFT's.
The main insight behind the program is the relation between these
Toeplitz operators and Hitchin's connection asymptotically in the quantum level
$k$. Here we have presented the initial application of this program, namely
the establishment of the asymptotic
faithfulness property for the quantum representations of the
mapping class groups. However this program can also be used to
study other asymptotic properties of these TQFT's. In particular
we have used them to establish that the quantum invariants for
closed 3-manifolds have asymptotic expansions in $k$. Topological
consequences of this are that certain classical topological
properties are determined by the quantum invariants, resulting in
interesting topological conclusions,
including very strong knot theoretical corollaries. Writeup of these
further developments of this program is in progress.

We remark that it is also an interesting
problem to understand how the Toeplitz operator constructions
used in this paper are related to the deformation quantization of
the moduli spaces described in \cite{AMR1} and \cite{AMR2}. In the
abelian case, the resulting Berezin-Toeplitz deformation quantization was
explicitly described in \cite{A2} and it turns out to be
equivalent to the one constructed in \cite{AMR2} in this case.

This paper is organized as follows. In section \ref{sec2} we give
the basic setup of applying geometric quantization to the moduli space to construct
the Verlinde bundle over \Teim space. In section \ref{sec3} we
review the construction of the connection in the
Verlinde bundle. We end that section by stating the properties
of the moduli space  and the Verlinde bundle we will use. We will
see that it is only a few elementary properties about the moduli
space, \Teim space and the general form of the connection in the
Verlinde bundle we really need.
 In section \ref{sec4} we review the
general
results about Toeplitz operators on smooth compact K\"{a}hler manifolds
we use in the following section
\ref{sec5}, where we prove that the Toeplitz operators for smooth
functions on the moduli space give asymptotically flat sections
of the endomorphism bundle of the Verlinde bundle. Finally, in
section \ref{sec6} we prove the asymptotic faithfulness
Theorem \ref{Main} above.

After the completion of this work Freedman and Walker, together
with Wang, have found a proof of the Asymptotic faithfulness
property for the SU(2)-BHMV-representations which uses
BHMV-technology. Their paper has already appeared \cite{FWW} (see
also \cite{M2} for a discussion). As alluded to before, we
are working jointly with K. Ueno to establish that these
representations are equivalent to our sequence $\rho_k^{2,0}$.
However, as long as this has not been established, our result is
logically independent of theirs.

For the SU(2)-BHMV-representations it is already known by the work of
Roberts \cite{Ro}, that they are irreducible for $k+2$ prime and that they have
infinite image by the work of Masbaum \cite{M}, except for a
few low values of $k$.

We would like to thank Nigel Hitchin, Bill Goldman and Gregor Masbaum for valuable
discussion. Further thanks is due to the Clay Mathematical
Institute for their financial support and to the Department of
Mathematics, University of California, Berkeley for their
hospitality, during the period, when this work was completed.

\section{The gauge theory construction of the Verlinde bundle}\label{sec2}

Let us now very briefly recall the construction of the Verlinde
bundle. Only the details needed in this paper will be given. We
refer e.g. to \cite{H} for further details.
As in the introduction we let $\Sigma$ be a closed oriented
surface of genus $g\geq 2$ and $p\in \Sigma$. Let $P$ be a
principal $SU(n)$-bundle over $\Sigma$. Clearly, all such $P$ are
trivializable. As above let $d\in \Z/n\Z \cong Z_{SU(n)}$. Throughout the rest of
this paper we will assume that $n$ and $d$ are coprime, although in the case $g=2$
we also allow $(n,d) = (2,0)$. Let $M$ be the moduli space of flat
$SU(n)$-connections in $P|_{\Sigma - p}$ with holonomy $d$ around
$p$. We can identify
\[M = \Hom_d({\tilde \pi}_1(\Sigma), SU(n) )/ SU(n).\]
Here ${\tilde \pi}_1(\Sigma)$ is the universal central extension
\[0 \ra \bZ \ra {\tilde \pi}_1(\Sigma) \ra \pi_1(\Sigma) \ra 1\]
as discussed in \cite{H} and in \cite{AB} and $\Hom_d$ means the
space of homomorphisms from ${\tilde \pi}_1(\Sigma)$
to $SU(n)$ which send the image of $1\in \bZ$ in
${\tilde \pi}_1(\Sigma)$ to $d$ (see \cite{H}). When
$n$ and $d$ are coprime, $M$ is a compact smooth manifold of dimension $m = (n^2
-1)(g-1)$. In general, when $n$ and $d$ are not coprime $M$ is not
smooth, except in the case where $g=2$, $n=2$ and $d=0$. In
this case $M$ is in fact diffeomorphic to  $ {\mathbb C}\P^3$.
There is a natural homomorphism from the mapping class group
to the outer automorphisms of ${\tilde \pi}_1(\Sigma)$, hence
$\Gamma$ acts on $M$.

We choose an invariant bilinear form $\{\cdot,\cdot\}$ on $\fg =
\mbox{Lie}(SU(n))$, normalized such that
$-\frac{1}{6}\{\vartheta\wedge [\vartheta\wedge\vartheta]\}$ is a
generator of the image of the integer cohomology in the real
cohomology in degree $3$ of $SU(n)$, where $\vartheta$ is the
$\fg$-valued Maurer-Cartan $1$-form on $SU(n)$.

This bilinear form induces a symplectic
form on $M$.  In fact
\[T_{[A]}M \cong H^1(\Sigma,d_A),\]
where $A$ is any flat connection in $P$ representing a point in $M$ and $d_A$ is the induced
covariant derivative in the associated adjoint bundle. Using this
identification, the symplectic form on $M$ is:
\[\omega(\varphi_1,\varphi_2) = \int_{\Sigma}\{\varphi_1\wedge\varphi_2\},\]
where $\varphi_i$ are $d_A$-closed $1$-forms on $\Sigma$ with
values in $\ad P$. See e.g. \cite{H} for further details on this.
The natural action of $\Gamma$ on $M$ is symplectic.

Let $\L$ be the Hermitian line bundle over $M$ and $\nabla$ the
compatible connection in $\L$ constructed by Freed in \cite{Fr}.
This is the content of Corollary 5.22, Proposition 5.24 and
equation (5.26) in \cite{Fr} (see also the work of Ramadas, Singer and Weitsman
\cite{RSW}). By Proposition 5.27 in \cite{Fr} we have that the
curvature of $\nabla$ is $\frac{\sqrt{-1}}{2\pi}\omega$. We will also use
the notation $\nabla$ for the induced connection in $\L^k$, where
$k$ is any integer.

By an almost identical construction, we can lift the action of
$\Gamma$ on $M$ to act on $\L$ such that the Hermitian connection
is preserved (See e.g. \cite{A1}). In fact, since $H^2(M, \bZ)\cong
\bZ$ and $H^1(M,\bZ) = 0$, it is clear that the action of $\Gamma$
leaves the isomorphism class of $(\L,\nabla)$ invariant, thus
alone from this one can conclude that a central extension of
$\Gamma$ acts on $(\L,\nabla)$ covering the $\Gamma$ action on
$M$. This is actually all we need in this paper, since we are only
interested in the projectivized action.

Let now $\sigma\in \T$ be a complex structure on $\Sigma$. Let us
review how $\sigma$ induces a complex structure on $M$ which is
compatible with the symplectic structure on this moduli space. The
complex structure $\s$ induces a $*$-operator on $1$-forms on $\Sigma$
and via Hodge theory we get that
\[H^1(\Sigma,d_A) \cong \ker(d_A + * d_A *).\]
On this kernel, consisting of the harmonic $1$-forms with values in
$\ad P$, the $*$-operator acts and its square is $-1$, hence we get an
almost complex structure on $M$ by letting $I= I_\sigma = -*$. It is a
classical result by Narasimhan and Seshadri (see \cite{NS1}),
that this actually makes $M$ a smooth
K{\"a}hler manifold, which as such, we denote $M_\sigma$. By
using the $(0,1)$ part of $\nabla$ in $\L$, we get an induced
holomorphic structure in the bundle $\L$. The resulting
holomorphic line bundle will be denoted
$\L_\s$. See also \cite{H} for further details on this.

From a more algebraic geometric point of view, we consider
the moduli space of S-equivalence classes of semi-stable bundles
of rank $n$ and determinant isomorphic to the line bundle ${\mathcal O}(d [p])$.
By using Mumford's Geometric
Invariant Theory, Narasimhan and Seshadri (see \cite{NS2}) showed
that this moduli space is a smooth complex algebraic projective variety
which is isomorphic as a
K{\"a}hler manifold to $M_\sigma$. Referring to \cite{DN} we
recall that

\begin{theorem}[Drezet \& Narasimhan]\label{DNTh1}
The Picard group of $M_\sigma$ is generated by the holomorphic line bundle
$\L_\s$ over $M_\sigma$
constructed above:
\[\Pic(M_\sigma) = \langle \L_\s \rangle. \]
\end{theorem}

\begin{definition}\label{Verlinde}
The {\em Verlinde bundle} $\V_k$ over \Teim space is by definition the
bundle whose fiber over $\sigma\in \T$ is $H^0(M_\sigma,\L_\s^k)$, where $k$ is a
positive integer.
\end{definition}

\section{The projectively flat connection}\label{sec3}
In this section we will review Axelrod, Della Pietra and Witten's
and Hitchin's construction of the projective flat connection over
\Teim space in the Verlinde bundle. We refer to \cite{H} and \cite{ADW}
for further details.

Let $\cH_k$ be the trivial $C^\infty(M,\L^k)$-bundle over $\T$
which contains $\V_k$, the Verlinde sub-bundle. If we have a
smooth one-parameter family of complex structures $\sigma_t$ on
$\Sigma$, then that induces a smooth one-parameter family of
complex structures on $M$, say $I_t$. In particular we get
$\sigma'_t\in T_{\sigma_t}(\T)$, which gives an $I'_t\in H^1(M_{\sigma_t},T)$
(here $T$ refers to the holomorphic tangent bundle of $M_{\sigma_t}$).

Suppose $s_t$ is a corresponding smooth one-parameter family in
$C^\infty (M,\L^k)$ such that $s_t \in H^0(M_{\sigma_t},\L_{\sigma_t}^k).$ By
differentiating the equation
\[(1+ \sqrt{-1} I_t) \nabla s_t = 0,\]
we see that
\[\sqrt{-1} I'_t\nabla s + (1+\sqrt{-1} I_t) \nabla s'_t = 0.\]
Hence, if we have an operator
\[u(v) : C^\infty(M,\L^k) \rightarrow C^\infty(M,\L^k)\]
for all real tangent vectors to \Teim space $v\in T(\T)$, varying smoothly with respect to $v$, and
satisfying
\[\sqrt{-1} I'_t \nabla^{1,0}s_t + \nabla^{0,1}u(\sigma'_t)(s_t) = 0,\]
for all smooth curves $\sigma_t$ in $\T$, then we get a connection induced in $\V_k$ by letting
\begin{equation}\Nabla_{v} = \Nabla^t_{v} - u(v),\label{Con}
\end{equation}
for all $v\in T(\T)$, where $\Nabla^t$ is the trivial connection
in $\cH_k$.

In order to specify the particular $u$ we are interested in, we
use
 the symplectic structure on $\omega\in \Omega^{1,1}(M_\sigma) $ to
define the tensor $G = G(v)\in \Omega^0(M_\s,S^2(T))$ by
\[v[I_\sigma] = G(v) \omega,\]
where $v[I_\sigma]$ means the derivative in the direction of $v\in  T_\s(\T)$ of the
complex structure $I_\sigma$ on $M$.
Following Hitchin, we give an explicit formula for $G$ in terms of
$v\in T_\s(\T)$:

The holomorphic tangent space to \Teim space at $\sigma\in \T$ is
given by
\[T^{1,0}_\sigma (\T) \cong H^1(\Sigma_\sigma, K^{-1}).\]
Furthermore, the holomorphic co-tangent space to the moduli space of
semi-stable bundles at the equivalence class of a stable bundle
$E$ is given by
\[T^*_{[E]}M_\sigma \cong H^0(\Sigma_\sigma, \End_0(E)\otimes K).\]
Thinking of $G(v)\in \Omega^0(M_\sigma,S^2(T))$ as a quadratic function on
$T^* = T^*M_\sigma$, we have that
\[G(v)(\alpha,\alpha) = \int_{\Sigma}\tr(\alpha^2)v^{(1,0)}\]
where $v^{(1,0)}$ is the image of $v$ under the projection onto $T^{1,0}(\T)$.
From this formula it is
clear that $G(v)\in H^0(M_\s,S^2(T))$ and that
$\Nabla$ agrees with $\Nabla^t$ along the anti-holomorphic
directions $T^{0,1}(\T)$. From Proposition (4.4) in \cite{H} we
have that this map $v\mapsto G(v)$ from $T_\s(\T)$ to $H^0(M_\s,S^2(T))$ is an
isomorphism.

The particular $u(v)$ we are interested in is $u_{G(v)}$, where
\begin{equation}
u_G(s) = \frac{1}{2(k + n)}(\Delta_G - 2  \nabla_{G \partial F} + \sqrt{-1} k f_G) s.
\label{Hitchcon}
\end{equation}

The leading order term $\Delta_G$ is the 2'nd order operator given by
\[
\begin{CD}
\Delta_G : \C^\infty(M,\L^k) @>{\nabla^{1,0}}>>
\C^\infty(M,T^*\otimes\L^k)
@>{G}>> \C^\infty(M,T\otimes\L^k)\\
 @>{\nabla^{1,0}\otimes 1 + 1\otimes\nabla^{1,0}}>>
 \C^\infty(M,T^*\otimes T\otimes \L^k) @>{\tr}>>
 \C^\infty(M,\L^k),
\end{CD}
\]
where we have used the Chern connection in $T$ on the K\"{a}hler
manifold $(M_\sigma,\omega)$.

The function $F=F_\s$ is the Ricci potential uniquely determined as the real
function with zero average over $M$, which satisfies the following
equation
\begin{equation}
\Ric_\sigma = 2 n \omega + 2 \sqrt{-1} \partial\bar{\partial} F_\sigma.
\label{riccipot}
\end{equation}
We usually drop the subscript $\s$ and think of $F$ as a smooth map from $\T$ to $\C^\infty(M)$.

The complex vector field $G\partial F \in \C^\infty(M_\s, T)$ is simply just the
contraction of $G$ with $\partial F \in \C^\infty(M_\s, T^*)$.

The
function $f_G\in C^\infty(M)$ is defined by
\[f_G = - \sqrt{-1} v[F],\]
where $v$ is determined by $G = G(v)$ and $v[F]$ means the derivative of
$F$ in
the direction of $v$. We refer to \cite{ADW} for this formula for
$f_G$.

\begin{theorem}[Axelrod, Della Pietra \& Witten; Hitchin]\label{Pflat}
The expression (\ref{Con}) above defines a connection in the bundle
$\V_k$, which induces a flat connection in $\P(\V_k)$.
\end{theorem}

Faltings has established this Theorem in the case where one
replaces $SU(n)$ with a general semi-simple Lie group (see
\cite{Fal}).

We remark about genus $2$, that \cite{ADW} covers this case, but
\cite{H} excludes this case, however, the work of Van Geemen and De
Jong \cite{vGdJ} extends Hitchin's approach to the genus $2$ case.

As discussed in the introduction, we see by Laszlo's Theorem that
this particular connection is the relevant one to study.

It will be essential for us to consider the induced flat
connection $\Nablae$ in the endomorphism bundle $\End(\V_k)$.
Suppose $\Phi$ is a section of $\End(\V_k)$. Then for all sections
$s$ of $\V_k$ and all $v\in T(\T)$ we have that
\[(\Nablae_v \Phi) (s) = \Nabla_v \Phi(s) - \Phi(\Nabla_v(s)).\]
Assume now that we have extended $\Phi$ to a section of $\Hom (\cH_k,\V_k)$ over $\T$.
Then
\begin{equation}\label{endocon}
\Nablae_v \Phi = \Nablaet_v \Phi - [\Phi, u(v)]
\end{equation}
where $\Nablaet$ is the trivial connection in the trivial
bundle $\End(\cH_k)$ over $\T$.

Let us end this section by summerazing the properties we use about
the moduli space in section \ref{sec5} to prove Theorem
\ref{Asympflat}, which intern implies Theorem \ref{Main}:

The moduli space $M$ is a finite dimensional
smooth compact manifold with a symplectic structure $\omega$, a
Hermitian line bundle $\L$ and a compatible connection $\nabla$,
whose curvature is $\frac{\sqrt{-1}}{2\pi}\omega$. Teichm\"{u}ller space
$\T$ is a smooth connected finite dimensional manifold, which
smoothly parametrizes K\"{a}hler structures $I_\sigma$,
$\sigma\in\T$ on $(M,\omega)$. For any positive integer $k$, we
have inside the  trivial bundle $\cH_k = \T \times \C^\infty(M,\L^k)$ the
finite dimensional subbundle $\V_k$, given by
\[\V_k(\sigma) = H^0(M_\sigma, \L_\sigma^k)\]
for $ \sigma \in \T$. We have a connection in $\V_k$ given by
\[\Nabla_v  = \Nabla^t_v  - u(v)\]
where $\Nabla^t_v$ is the trivial connection in $\cH_k$
and $u(v)$ is the second order differential
operator $u_{G(v)}$ given in (\ref{Hitchcon}). - All we will need
about the operator $\Delta_G - 2 \nabla_{G\partial F}$, is that
there is some finite set of vector
fields  $X_r(v), Y_r(v), Z(v) \in
\C^\infty(M_\sigma,T)$, $r = 1, \ldots, R$ (where $v\in T_\sigma(\T)$) all varying
smoothly\footnote{This makes sense when we consider the
holomorphic tangent bundle $T$ of $M_\sigma$ inside the
complexified real tangent
bundle $TM \otimes {\mathbb C}$ of $M$.} with
$v\in T(\T)$, such that
\begin{equation}
\Delta_{G(v)} -2 \nabla_{G(v)\partial F} = \sum_{r=1}^R \nabla_{X_r(v)}
\nabla_{Y_r(v)} + \nabla_{Z(v)}.\label{Opform}
\end{equation}
This follows immediately from the definition of $\Delta_{G(v)}$.
From this we have the expression
\begin{equation}
u(v) = \frac1{2(k+n)}\left( \sum_{r=1}^R \nabla_{X_r(v)}
\nabla_{Y_r(v)} + \nabla_{Z(v)} + n v[F]\right)  - \frac12 v[F].
\label{Hitchcons}
\end{equation}
All we need to use about $F : \T \ra C^\infty(M)$ is that it is a
smooth function, such that $F_\sigma$ is real valued
on $M$ for all $\sigma\in \T$.

\section{Toeplitz operators on compact K{\"a}hler manifolds}\label{sec4}

In this section $(N^{2m},\omega)$ will denote a compact K{\"a}hler
manifold on which we have a holomorphic line bundle $L$ admitting
a compatible Hermitian connection whose curvature is $\frac{\sqrt{-1}}{2\pi}\omega$.
On $\C^\infty(N,L^k)$ we have the $L_2$-inner product:
\[\langle s_1, s_2 \rangle = \frac{1}{m!}\int_N (s_1,s_2) \omega^m\]
where $s_1, s_2 \in \C^\infty(N,L^k)$ and $(\cdot,\cdot)$ is the fiberwise
Hermitian structure in $L^k$. This $L_2$-inner product gives the
 orthogonal projection $\pi : \C^\infty(N,L^k) \ra H^0(N,L^k)$.
For each $f\in \C^\infty(N)$ consider
the associated {\em Toeplitz operator} $T_f^{(k)}$ given as the
composition of the multiplication operator $M_f : H^0(N,L^k) \ra
\C^\infty(N,L^k)$ with the orthogonal projection $\pi :
\C^\infty(N,L^k) \ra H^0(N,L^k)$, so that
\[T_f^{(k)}(s) = \pi(f s).\]
Since the multiplication operator is a zero order differential
operator, $T_f^{(k)}$ is a zero-order Toeplitz operator.
Sometimes we will suppress the
superscript $(k)$ and just write $T_f = T_f^{(k)}$.

Let us here
give an explicit formula for $\pi$:  Let $h_{ij} = \langle s_i,
s_j\rangle$, where $s_i$ is a basis of $H^0(N,L^k)$. Let $h^{-1}_{ij}$ be the
inverse matrix of $h_{ij}$. Then
\begin{equation}\label{projf}
\pi(s) = \sum_{i,j}\langle s, s_i\rangle h^{-1}_{ij}s_j.
\end{equation}
This formula will be useful when we have to compute the
derivative of $\pi$ along a one-parameter curve of complex
structures on the moduli space.

Suppose we have a smooth section $X\in C^\infty(N, TN)$ of the holomorphic tangent
bundle of $N$. We then claim that the operator $\pi \nabla_X$
is a zero-order Toeplitz operator. Suppose $s_1\in C^\infty(N,L^k)$ and $s_2 \in
H^0(N,L^k)$, then we have that
\[X(s_1,s_2) = (\nabla_X s_1, s_2).\]
Now, calculating the Lie derivative along $X$ of $(s_1,s_2)\omega^m$ and using
 the above,
one obtains after integration that
\[\langle \nabla_X s_1, s_2 \rangle = - \langle \Lambda d(i_X\omega) s_1, s_2 \rangle,\]
where $\Lambda$ denotes contraction with $\omega$. Thus
\begin{equation}\pi \nabla_X = T_{f_X}^{(k)},\label{1to0order}
\end{equation}
as operators from $C^\infty(N,L^k)$ to $H^0(N,L^k)$, where $f_X = -\Lambda d(i_X\omega)$.
Iterating this, we find for all $X_1,X_2 \in C^\infty(TN)$ that
\begin{equation}\pi \nabla_{X_1}\nabla_{X_2} = T^{(k)}_{f_{X_2}f_{X_1} - X_2(f_{X_1})}\label{2to0order}
\end{equation}
again as operators from $C^\infty(N,L^k)$ to $H^0(N,L^k)$.

For $X\in C^\infty(N, TN)$, the complex conjugate vector field
$\bar X \in \C^\infty (N, \bar T N)$ is a section of the antiholomorphic tangent bundle,
and for $s_1, s_2 \in C^\infty(N,L^k)$, we have that
\[\bar X (s_1,s_2) = (\nabla_{\bar X} s_1, s_2) + ( s_1, \nabla_X s_2).\]
Computing the Lie derivative along $\bar X$ of
$(s_1,s_2)\omega^m$ and integrating, we get that
\[\langle \nabla_{\bar X} s_1, s_2 \rangle + \langle (\nabla_{X})^* s_1, s_2 \rangle
 = \langle \Lambda d(i_{\bar X}\omega) s_1, s_2 \rangle.\]
Hence we see that
\[(\nabla_{X})^* = - \left( \nabla_{\bar X} + f_{\bar X} \right)\]
as operators on $C^\infty(N,L^k)$. In particular, we see that
\begin{equation}\label{1to0order*}
\pi (\nabla_{X})^* \pi = - T_{f_{\bar X}}|_{H^0(N,L^k)} : H^0(N,L^k) \ra H^0(N,L^k).
\end{equation}
For two smooth sections $X_1, X_2$ of the holomorphic tangent bundle
$TN$ and a smooth function $h\in C^\infty(N)$, we deduce from the
formula for $(\nabla_{X})^*$ that
\begin{eqnarray}\label{2to0order*}
\pi (\nabla_{{X}_1})^*(\nabla_{{X}_2})^* h\pi & = & \pi \bar X_1 \bar X_2(h)\pi +\\
& & \quad
\pi f_{{\bar X}_1} \bar X_2(h)\pi +\pi f_{{\bar X}_2}\bar X_1(h) \pi +
 \nonumber \\
& & \quad \pi\bar X_1(f_{{\bar X}_2}) h \pi + \pi f_{{\bar X}_1} f_{{\bar
X}_2}h
\pi\nonumber
\end{eqnarray}
as operators on $H^0(N,L^k)$.

We need the following Theorems on Toeplitz operators. The first is
due to Bordemann, Meinrenken and Schlichenmaier (see \cite{BMS}).
The $L_2$-inner product on $\C^\infty(N,L^k)$ induces an inner product on $H^0(N,L^k)$, which in
turn induces the operator norm $\|\cdot\|$ on $\End(H^0(N,L^k))$.

\begin{theorem}[Bordemann, Meinrenken and Schlichenmaier]\label{BMS1}
For any $f\in \C^\infty(N)$ we have that
\[\lim_{k\ra \infty}\|T_{f}^{(k)}\| = \sup_{x\in N}|f(x)|.\]
\end{theorem}

Since the association of the sequence of Toeplitz operators
$T^k_f$, $k\in \Z_+$ is linear in $f$, we see from this Theorem,
that this association is faithful.

\begin{theorem}[Schlichenmaier]\label{S}
For any pair of smooth functions $f_1, f_2\in \C^\infty(N)$, we
have an asymptotic expansion
\[T_{f_1}^{(k)}T_{f_2}^{(k)} \sim \sum_{l=0}^\infty T_{c_l(f_1,f_2)}^{(k)} k^{-l},\]
where $c_l(f_1,f_2) \in C^\infty(N)$ are uniquely determined since
$\sim$ means the following: For all $L\in \Z_+$ we have that
\begin{equation}
\|T_{f_1}^{(k)}T_{f_2}^{(k)} - \sum_{l=0}^L T_{c_l(f_1,f_2)}^{(k)} k^{-l}\| =
O(k^{-(L+1)}).\label{normasympToep}
\end{equation}
Moreover, $c_0(f_1,f_2) = f_1f_2$.
\end{theorem}

This Theorem
was proved in \cite{Sch} and is published in \cite{Sch1} and \cite{Sch2}, where it
is also proved that the formal
generating series for the $c_l(f_1,f_2)$'s gives a formal
deformation quantization of the Poisson structure on $N$ induced from
$\omega$. By examining the proof in \cite{Sch} (or in \cite{Sch1} and \cite{Sch2})
of this Theorem, one observes
that for continuous families of functions, the
estimates in Theorem \ref{S} are uniform over compact parameter spaces.

\section{Toeplitz operators on moduli space and the projective flat connection}\label{sec5}

Let $f\in \C^\infty(M)$ be a smooth function on the moduli space.
We consider $T_f^{(k)}$ as a section
of the endomorphism bundle $\End(\V_k)$. The flat connection
$\Nabla$ in the projective bundle $\P(\V_k)$ induces the flat
connection $\Nablae$ in the endomorphism bundle $\End(\V_k)$ as described in section \ref{sec3}. We
shall now establish that the sections $T_f^{(k)}$ are in a certain
sense asymptotically flat by proving the following Theorem.

\begin{theorem}\label{Asympflat}
Let $\sigma_0$ and $\sigma_1$ be two points in \Teim space and
$P_{\sigma_0,\sigma_1}$ be the parallel transport in the flat
bundle $\End(\V_k)$ from $\sigma_0$ to $\sigma_1$. Then
\[\|P_{\sigma_0,\sigma_1}T_{f,\sigma_0}^{(k)} - T_{f,\sigma_1}^{(k)}\| = O(k^{-1}),\]
where $\|\cdot\|$ is the operator norm on $H^0(M_{\s_1},\L_{\s_1}^k)$.
\end{theorem}

In the proof of this Theorem we will make use of the following Hermitian structure on $\cH_k$
\begin{equation}
\langle s_1,s_2\rangle_F = \frac{1}{m!}\int_M (s_1,s_2) e^{- F}\omega^m, \label{asymphermform}
\end{equation}
where we recall that $F=F_\s$ is the Ricci potential, which is a
real smooth function on $M_\s$ for each $\s\in \T$ determined by equation
(\ref{riccipot}).
In Lemma \ref{equinormL2} below we will see that $\langle
\cdot,\cdot\rangle_F$ is uniform equivalent to the constant $L_2$-Hermitian
structure on $\cH_k$,
when both are restricted to $\V_k$ over any compact subset of $\T$.
The constant $L_2$-Hermitian structure on $\cH_k$ is
not asymptotically flat with respect to $\Nabla$, but
Proposition \ref{asympherm} below shows that the Hermitian structure
$\langle \cdot,\cdot\rangle_F$
restricted to $\V_k$ is asymptotically flat with respect to $\Nabla$. - It therefore
induces a Hermitian structure on $\V_k^*\otimes \V_k$ (which we also
denote $\langle\cdot,\cdot\rangle_F$), which is asymptotic flat with respect to $\Nablae$.

This suggests one considers the smooth function $t\mapsto |P_{\sigma_0,\sigma_t}
T_{f,\sigma_0}^{(k)} - T_{f,\sigma_t}^{(k)}|_F$ and establishes an
$O(k^{-1})$ estimate for its derivative uniformly over $J$, since
by Lemma \ref{equinormL2} and the first inequality in (\ref{comp})
below,
$O(k^{-1})$ control on $|P_{\sigma_0,\sigma_1}
T_{f,\sigma_0}^{(k)} - T_{f,\sigma_1}^{(k)}|_F$ implies Theorem \ref{Asympflat}. An
$O(k^{-1})$ estimate on
$|\Nablae_{\sigma'_t}T_{f,\sigma_t}^{(k)}|_F$ uniformly over $J$
would imply this estimate, however we are only able to establish
that $\|\Nablae_{\sigma'_t}T_{f,\sigma_t}^{(k)}\|$ is $O(k^{-1})$
uniformly over the interval, which is considerably weaker because
of the $\frac{m}{2}$'th
power of $k$ in the second inequality in (\ref{comp}) below.

We shall therefore perturb the function $f$ by
adding on sufficiently many terms of the form $(h_l)_t k^{-l}$
($l=1, ..., r>m/2$), where $(h_l)_t \in C^\infty(M)$, $t\in J$, so as to obtain a smooth
one-parameter family of functions $(f_r)_t =
f + \sum_{l=1}^r (h_l)_t k^{-l}$. The
$(h_l)_t$'s are determined inductively in $l$ such that
$\|\Nabla_{\sigma'_t}^e T_{(f_r)_t,\sigma_t}^{(k)}\|$
is $O(k^{-r-1})$ uniformly over $J$. This is the content of Proposition
\ref{Mainnorm} below. That estimate on the covariant derivative of
the Toeplitz operator $T^{(k)}_{(f_r)_t,\sigma_t}$
will allow us to prove that $\|P_{\sigma_0,\sigma_1}T_{(f_r)_0,\sigma_0}^{(k)} -
T_{(f_r)_1,\sigma_1}^{(k)}\|$ is $ O(k^{-1})$ by analyzing the derivative of
\[t \mapsto |P_{\sigma_0,\sigma_t}T^{(k)}_{(f_r)_0,\sigma_0} - T^{(k)}_{(f_r)_t,\sigma_t}|^2_F.\]
Since we can arrange that
$(f_r)_0 = f$ and since $\| T_{(f_r)_1,\sigma_1}^{(k)}-
T_{f,\sigma_1}^{(k)}\|$ is $O(k^{-1})$, this will allow us to
prove Theorem \ref{Asympflat}.

First however, we need to establish a useful formula for the derivative of the
orthogonal projection $\pi$ along the curve $\sigma_t$. To this end, consider a basis of
covariant constant sections\footnote{We will from now on mostly
suppress the subscript $t$ on various
quantities defined along the curve $\sigma_t$.} $s_i = (s_i)_t$, $i=1,\ldots, \rank \V_k$, of
$\V_k$ over the curve $\sigma_t$:
\[s'_i = u_G(s_i), \qquad i=1, \ldots , \rank \V_k.\]
Recall formula (\ref{projf}) for the projection $\pi : \C^\infty(M,\L^k) \ra
H^0(M_{\s_t},\L_{\s_t}^k)$ and compute the derivative along $\s_t$: For
any fixed $s\in \C^\infty(M,\L^k)$, we have that
\begin{eqnarray*}
\pi'(s) & = & \sum_{i,j}\langle s, s'_i\rangle h^{-1}_{ij}s_j\\
&& + \sum_{i,j}\langle s, s_i\rangle (h^{-1}_{ij})'s_j\\
&& + \sum_{i,j}\langle s, s_i\rangle h^{-1}_{ij}s'_j.
\end{eqnarray*}
An easy computation gives that
\[(h^{-1}_{ij})' = -\sum_{l,r}h^{-1}_{il}(\langle s'_l, s_r\rangle +
\langle s_l, s'_r\rangle) h^{-1}_{rj},\]
so
\begin{eqnarray*}
\pi\pi'(s) & = & \sum_{i,j}\langle u_G^* s, s_i\rangle
h^{-1}_{ij}s_j\\
&& - \sum_{i,l,m,j}\langle s, s_i\rangle h^{-1}_{il}\langle s_l, s'_m\rangle
h^{-1}_{mj}s_j\\
& = & \pi u_G^*(s) - \sum_{m,j}\langle \pi s,  s'_m\rangle h^{-1}_{mj}s_j\\
& = & \pi  u_G^*(s) - \pi  u_G^* \pi(s).
\end{eqnarray*}
Hence we conclude that
\begin{equation}
\pi \pi' = \pi u_G^* - \pi u_G^*\pi. \label{derivepi}
\end{equation}

Having derived the formula for the derivative of $\pi$, we now
proceed to construct the needed perturbation of $f$. Let $r$ be a non-negative integer and let
\[(f_r)_t = \sum_{i=0}^r (h_i)_t k^{-i},\]
where the $(h_i)_t$, $t\in J$, for now are arbitrary
smooth one-parameter families of smooth functions
on $M$, however we will fix $(h_0)_t =
f$ for all $t\in J$.
We have that $T_{f_r}^{(k)}$ is a section of $\End(\V_k)$ over the curve $\sigma_t$.
According to formula (\ref{endocon}), we have that
\[\Nabla_{\sigma'}^e(T_{f_r})= (T_{f_r})' - [u_G,T_{f_r}]
= \pi f'_r + \pi' f_r  - [u_G,\pi f_r].\]
Since $\Nabla_{\sigma'}^e(T_{f_r})$ is a section of $\End(\V_k)$,
we have of course that
\[\pi \Nablae_{\sigma'}(T_{f_r}) \pi = \Nabla_{\sigma'}^e(T_{f_r})
\pi :H^0(M_{\sigma_t},\L_{\s_t}^k) \ra H^0(M_{\sigma_t},\L_{\s_t}^k).\]

\begin{proposition}\label{Mainnorm} Given $f$ and a non-negative integer $r$,
there exists unique smooth one-parameter families
of functions $(h_i)_t \in  C^\infty(M)$, $i=1,\ldots, r$ and $t\in J$ such that
\begin{equation}
\sup_{t\in J}\|\Nablae_{\sigma'}(T_{f_r})\| = O(k^{-r-1}),\label{covconstm}
\end{equation}
and $(h_i)_0 = 0$, $i=1,\ldots, r$.
\end{proposition}

\proof For any choice of $h_i$'s, $i=1,\ldots, r$, we compute using (\ref{derivepi}) that
\[
\begin{split}
\pi \Nabla_{\sigma'}^e(T_{f_r}) \pi & = \sum_{i=1}^r \pi h'_i \pi k^{-i}\\
            & \quad + \sum_{i=0}^r (\pi u_G^* h_i\pi - \pi u_G^* \pi h_i\pi) k^{-i}\\
            & \quad - \sum_{i=0}^r (\pi u_G \pi h_i\pi - \pi h_i u_G\pi) k^{-i}.
\end{split}\]

Using (\ref{Opform}) together with (\ref{1to0order}) and (\ref{2to0order}) we
define the
function $H_G\in \C^\infty(M)$ independent of $k$, as follows
\[\pi (\Delta_G - 2 \nabla_{G \partial F}) = \pi H_G.\]
We also define $H_i\in C^\infty(M)$ which depend on $h_i$ and $G$ but are independent of
 $k$, by
\[\pi h_i(\Delta_G - 2 \nabla_{G \partial F}) = \pi H_i.\]
Using (\ref{1to0order*}) and (\ref{2to0order*}) we
further define the
function $H^*_G\in \C^\infty(M)$ independent of $k$, as follows
\[\pi (\Delta_G^* - 2 (\nabla_{G \partial F})^*)\pi = \pi H_G^*\pi.\]
And we define $H_i^*\in C^\infty(M)$ which depend on $h_i$ and $G$ but are independent of
 $k$, by
\[\pi (\Delta_G^* - 2 (\nabla_{G \partial F})^*)h_i\pi = \pi H_i^* \pi.\]
Then we have that
\[
\begin{split}
\pi \Nabla_{\sigma'}^e(T_{f_r}) \pi & = \sum_{i=1}^r \pi h'_i \pi k^{-i}\\
            & \quad - \sum_{i=0}^r \frac1{2(k+n)}(\pi H_G^* \pi h_i\pi
            - \pi H_i^*\pi) k^{-i}\\
            & \quad - \sum_{i=0}^r \frac1{2(k+n)}(\pi H_G \pi h_i\pi - \pi H_i\pi)
            k^{-i}\\
            & \quad + \sum_{i=0}^r \frac{k}{2(k+n)}\sqrt{-1}( \pi \overline{f}_G \pi h_i\pi - \pi \overline{f}_G h_i\pi) k^{-i}\\
            & \quad - \sum_{i=0}^r \frac{k}{2(k+n)}\sqrt{-1}(\pi f_G \pi h_i\pi - \pi h_i f_G\pi)
            k^{-i}.
\end{split}\]

Because of (\ref{normasympToep}) and Theorem \ref{BMS1}, we see
for any $r$,
that condition (\ref{covconstm}) is equivalent to

\[
\begin{split}
\sup_{t\in J}\| & \sum_{i=1}^r \pi h'_i \pi k^{-i}\\
            & \quad - \frac1{2k}\sum_{i=0}^m \sum_{j=0}^r (-1)^j n^j
           (\sum_{l=0}^r \pi c_l(H_G^*, h_i)\pi k^{-l}
            - \pi H_i^*\pi) k^{-i-j}\\
            & \quad -  \frac1{2k}\sum_{i=0}^r \sum_{j=0}^r (-1)^j n^j
           (\sum_{l=0}^r \pi c_l(H_G, h_i)\pi k^{-l} - \pi H_i\pi)
            k^{-i-j}\\
            & \quad +  \frac{\sqrt{-1}}{2k}\sum_{i=0}^r \sum_{j=0}^r
           \sum_{l=1}^r (-1)^j n^j \pi c_l(\overline{f}_G, h_i)\pi k^{-l-i-j+1} \\
            & \quad - \frac{\sqrt{-1}}{2k}\sum_{i=0}^r \sum_{j=0}^r
           \sum_{l=1}^r (-1)^j n^j \pi c_l(f_G, h_i)\pi k^{-l-i-j+1} \| = O(k^{-r-1}).
\end{split}\]

Note that the expression inside the norm is a polynomial in
$k^{-1}$, say $\sum_{i=1}^{N_r} \pi (C_i)_t\pi k^{-i}$, where each coefficient
$\pi (C_i)_t \pi$
is the Toeplitz operator associated to a smooth one-parameter family
$(C_i)_t\in C^\infty(M)$, $t\in J$. We also note that the
functions $C_i$, $i=1, \ldots, r$ do not depend on $r$ and in
the same range for $i$ they are all of the following form
\[C_i = h'_i - D_i\]
where $(D_i)_t\in C^\infty(M)$, $t\in J$, is a linear combination
of
$$H_j^*, H_j, c_l(H_G^*,h_j), c_l(H_G,h_j), c_l(\bar{f}_G,h_j),
c_l(f_G,h_j) \in C^\infty(J, C^\infty(M)),$$
for $j = 0,\ldots, i-1$  and $l=0, \ldots, i$, e.g. for
$i=1$ we get
\[D_{1} = - \frac12 (H_G^* h_0 - H_0^* - H_G h_0  - H_0 - \sqrt{-1}
c_1(\bar f_G, h_0) + \sqrt{-1} c_1(f_G, h_0)). \]

Now observe by Theorem \ref{BMS1} that condition (\ref{covconstm})
is equivalent to $C_i = 0$ for $i=1, \ldots r$.
From this it follows that we can inductively uniquely
determine the functions $h_i$. First of all, since the
coefficient of $k^0$ is zero, we see that (\ref{covconstm}) holds
for $r=0$.  Now assume that we have determined $h_1, \ldots
h_{r-1}$ uniquely such that $(h_i)_0= 0$ and
\[ h'_i = D_i,\]
for $i = 1, \ldots , r-1$.
By the above we then observe that condition (\ref{covconstm}) holds if
and only if $C_r = 0$, which means if and only if
\[h'_r = D_r.\]
But since we require $(h_r)_0 = 0$, we must have
that
\[(h_r)_t = \int_0^t (D_r)_s ds.\]
\eproof

\begin{lemma}\label{equinormL2}
The Hermitian structure on $\cH_k$
$$\langle s_1,s_2\rangle_F = \frac{1}{m!}\int_M (s_1,s_2)
e^{- F}\omega^m $$
and the constant $L_2$-Hermitian structure on $\cH_k$
$$\langle s_1,s_2 \rangle = \frac{1}{m!}\int_M (s_1,s_2) \omega^m$$
are equivalent uniformly in $k$ when restricted to $\V_k$
over any compact subset $K$ of $\T$.
\end{lemma}

\proof

We clearly have that
\[|s|^2_F  \leq \|T^{(k)}_{e^{- F }}\|\ |s|^2, \]
so by Theorem \ref{BMS1} we see that there exists a constant $C$
(depending on K)
such that
\[|s|_F  \leq C |s|, \]
for all $k$. Conversely, we have that
\[\begin{split}
|s|^2 & = \langle \pi e^{\frac12 F}e^{-\frac12 F}\pi s, s\rangle\\
&\leq |\langle (\pi e^{\frac12 F} e^{-\frac12 F }\pi - \pi e^{\frac12 F}
\pi e^{-\frac12 F }\pi)s,s\rangle | \\
&\quad + |\langle  \pi e^{\frac12 F}
\pi e^{-\frac12 F }\pi s,s\rangle | \\
& \leq \| \pi e^{\frac12 F} e^{-\frac12F }\pi - \pi e^{\frac12 F}
\pi e^{-\frac12 F }\pi \|\ |s|^2\\
&\quad +  \| \pi e^{\frac12 F}\pi \|\ |\pi e^{-\frac12 F}\pi s |\ |s|.
\end{split}\]
By Theorem \ref{BMS1} and \ref{S} we see there exists constants
$C'$ and $C''$ (again depending on K) such that
\[|s| \leq \frac{C'}k |s| + C'' |s|_F.\]
But then we have for all sufficiently large $k$ that
\[|s| \leq 2 C'' |s|_F.\]
Hence we have established the claimed equivalence.
\eproof

Along any smooth one-parameter family of complex structures
$\sigma_t$, we have that
\[\frac{d}{d t} \langle s_1,s_2\rangle_F =
\langle \Nabla^t_{\sigma'_t} s_1,s_2\rangle_F + \langle
s_1,\Nabla^t_{\sigma'_t} s_2\rangle_F - \langle
\frac{\partial F}{\partial t} s_1,s_2\rangle_F.\] So,
if we let
\[E(s) = \frac{d}{d t} |s|^2_F -
\langle \Nabla_{\sigma'_t}s,s\rangle_F -
\langle s,\Nabla_{\sigma'_t}s\rangle_F \]
recalling that $\sqrt{-1}f_G = \frac{\partial F}{\partial t}$, $\Nabla_v = \Nabla^t_v - u(v)$
and formula (\ref{Hitchcons}),
 we have for all sections $s$ of $\V_k$ that
\begin{eqnarray*}
E(s) & = & \frac{1}{2(k+n)}
\left(\langle \pi e^{- F}( \Delta_G - 2 \nabla_{G \partial F} -\sqrt{-1} n f_G) s, s\rangle\right. \\
& & \left. \phantom{hjdfhkdfdffh} + \langle  s,
\pi e^{-F}(\Delta_G - 2 \nabla_{G \partial F} -\sqrt{-1} n f_G) s\rangle\right).
\end{eqnarray*}
Hence by combining Theorem \ref{BMS1}, (\ref{Opform}),  (\ref{1to0order}) and (\ref{2to0order})
we have proved that
\begin{proposition}\label{asympherm}
The Hermitian structure (\ref{asymphermform}) is asymptotically flat
with respect to the connections $\Nabla$, i.e. for any compact subset $K$ of
$\T$, there exists a
constant $C$ such that for all sections $s$ of $\V_k$ over $K$, we
have that
\[|E(s)| \leq \frac{C}{k+n} |s|_F^2\]
over $K$.
\end{proposition}

We note that this proposition implies the same proposition for
sections of $\End(\V_k)$ with respect to the induced
Hermitian structure on $\End(\V_k)= \V_k^*\otimes \V_k$, which we also
denote $\langle\cdot,\cdot\rangle_F$.
We denote the analogous quantity of $E$ for the endomorphism bundle
by $E_e$.

\proof[Proof of Theorem \ref{Asympflat}]
Let $\sigma_t$, $t\in J$  be a smooth one-parameter family of complex structures
such that $\sigma_t$ is a curve in
$\T$ between the two points in question.
By Lemma \ref{equinormL2} the Hermitian structures
$\langle\cdot,\cdot\rangle$ and $\langle\cdot,\cdot\rangle_F$ on $\V_k$ are
equivalent uniformly in $k$ over compact subsets of $\T$. It follows there exists constants
$C_1$ and $C_2$, such that we have the following inequalities over
the image of $\sigma_t$ in $\T$ for the operator norm $\|\cdot\|$
and the norm $| \cdot |_F$ on $\End(\V_k)$:
\begin{equation}\|\cdot\| \leq C_1 | \cdot |_F \leq
C_2 \sqrt{P_{g,n}(k)}\|\cdot\|,\label{comp}\end{equation}
where $P_{g,n}(k)$ is the rank of $\V_k$ given by the Verlinde formula.
By the Riemann-Roch Theorem this is a polynomial in $k$ of degree $m$.

Because of these inequalities, we choose an integer $r $ bigger
than $m/2$ and let $f_r$ be as provided by Proposition
\ref{Mainnorm} and we define $n_k : J \ra [0,\infty)$ by
\[n_k(t) = |\Theta_k(t)|_F^2\]
where
\[\Theta_k(t) : \V_{k, \sigma_t} \ra \V_{k, \sigma_t}\]
is given by
\[ \Theta_k(t) = P_{\sigma_0,\sigma_t}T^{(k)}_{(f_r)_0,\sigma_0} - T^{(k)}_{(f_r)_t,\sigma_t}.\]
The functions $n_k$ are differentiable in $t$ and we compute that
\[
\begin{split}
\frac{d n_k}{d t} & = \langle
\Nabla_{\sigma'_t}^e(\Theta_k(t)),\Theta_k(t)\rangle_F
+ \langle
\Theta_k(t),
\Nabla_{\sigma'_t}^e(\Theta_k(t))\rangle_F + E_e(\Theta_k(t))\\
&=- \langle \Nabla_{\sigma'_t}^eT^{(k)}_{(f_r)_t,\sigma_t},
\Theta_k(t)\rangle_F - \langle
\Theta_k(t),
\Nabla_{\sigma'_t}^e T^{(k)}_{(f_r)_t,\sigma_t}\rangle_F + E_e(\Theta_k(t)).
\end{split}
\]
Using the above, we get the following estimate
\[
\begin{split}
|\frac{d n_k}{d t}| &\leq 2 |\Nabla_{\sigma'_t}^e
T^{(k)}_{(f_r)_t,\sigma_t}|_F |\Theta_k(t)|_F + |E_e(\Theta_k(t))|\\
& \leq 2 C_2 \sqrt{P_{g,n}(k)} \|\Nabla_{\sigma'_t}^e
T^{(k)}_{(f_r)_t,\sigma_t}\| n_k^{1/2} + |E_e(\Theta_k(t))|.\end{split}
\]
Consequently we can apply Proposition \ref{Mainnorm} and \ref{asympherm} to obtain that there
exists a constant $C$ such that
\[
|\frac{d n_k}{d t}| \leq \frac{C}{k} ( n_k^{1/2} + n_k).
\]
This estimate implies that
\[n_k(t) \leq (\exp(\frac{C t}{2k}) - 1 )^2.\]
But by (\ref{comp}) we get that
\[\|P_{\sigma_0,\sigma_1}T_{(f_r)_0,\sigma_0}^{(k)} - T_{(f_r)_1,\sigma_1}^{(k)}\|
= \|\Theta_k(1)\|\leq  C_1 n_k(1)^{1/2}.\]
The Theorem then follows from these two estimates, since
$(f_r)_0 = f$ and
$$\|T_{(f_r)_1,\sigma_1}^{(k)}- T_{f,\sigma_1}^{(k)} \| =
O(k^{-1}).$$
\eproof

\section{Asymptotic faithfulness}\label{sec6}

Recall that the flat connection in the bundle $\P(\V_k)$ gives the
projective representation of the mapping class group
\[\rho_k : \Gamma \ra \Aut(\P(V_k)),\]
where $\P(V_k) = $ covariant constant sections of $\P(\V_k)$ over \Teim space.

\begin{theorem}\label{Main1}
For any $\phi\in \Gamma$, we have that
\[\phi \in \bigcap_{k=1}^\infty \ker \rho_k\]
if and only if $\phi$ induces the identity on $M$.
\end{theorem}

\proof
Suppose we have a $\phi \in \Gamma$. Then $\phi$ induces a symplectomorphism of $M$
which we also just denote $\phi$ and we get the following commutative
diagram for any $f\in \C^\infty(M)$
\begin{equation*}
\begin{CD}
H^0(M_\sigma,\L_{\s}^k) @>\phi^*>> H^0(M_{\phi(\sigma)},\L_{\phi(\sigma)}^k)
@> P_{\phi(\sigma),\sigma}>> H^0(M_\sigma,\L_{\s}^k)\\
@V T^{(k)}_{f,\sigma} VV @V T^{(k)}_{f \circ \phi,\phi(\sigma)} VV
@VV{P_{\phi(\sigma),\sigma}T^{(k)}_{f \circ \phi,\phi(\sigma)}}V\\
H^0(M_\sigma,\L_{\s}^k) @>\phi^*>> H^0(M_{\phi(\sigma)},\L_{\phi(\sigma)}^k)
@> P_{\phi(\sigma),\sigma}>> H^0(M_\sigma,\L_{\s}^k),
\end{CD}
\end{equation*}
where $P_{\phi(\sigma),\sigma} : H^0(M_{\phi(\sigma)},\L_{\phi(\sigma)}^k) \ra H^0(M_{\sigma},\L_{\s}^k)$
on the horizontal arrows refer to parallel transport in the Verlinde
bundle itself, whereas  $P_{\phi(\sigma),\sigma}$  refers to the parallel transport in
the endomorphism bundle $\End(\V_k)$ in the last vertical arrow.
Suppose now $\phi \in \bigcap_{k=1}^\infty \ker
\rho_k$, then $P_{\phi(\sigma),\sigma} \circ \phi^* = \rho_k(\phi) \in
\bC \id$ and we get that $T^{(k)}_{f,\sigma} = P_{\phi(\sigma),\sigma}
T^{(k)}_{f\circ \phi,\phi(\sigma)}$. By Theorem \ref{Asympflat} we get
that
\[
\begin{split}
\lim_{k\ra \infty}\|T_{f - f \circ \phi,\sigma}^{(k)}\| &=
\lim_{k\ra \infty}\|T_{f,\sigma}^{(k)} - T_{f\circ \phi,\sigma}^{(k)}\|\\
& =  \lim_{k\ra \infty}\| P_{\phi(\sigma),\sigma}
T^{(k)}_{f\circ \phi,\phi(\sigma)} - T_{f\circ \phi,\sigma}^{(k)} \| = 0.
\end{split}\]
By Bordemann, Meinrenken and Schlichenmaier's Theorem \ref{BMS1},
we must have that $f = f \circ \phi$. Since this holds for any
$f\in \C^\infty(M)$, we must have that $\phi$ acts by the identity
on $M$. \eproof

\proof[Proof of Theorem \ref{Main}]
Our main Theorem \ref{Main} now follows directly from Theorem \ref{Main1},
since it is known that the only
element of $\Gamma$, which acts by the identity on the moduli space
$M$ is the identity, if $g>2$. If $g=2$, $n=2$ and $d$ is even,
it is contained in the sub-group generated
by the hyper-elliptic involution, else it is also the identity.

A way to see this using the moduli space of
flat $SL(n,\bC)$-connections goes as follows:
We have that $M$ is a component of the
real slice in $\M$, the moduli space of flat $SL(n,\C)$-connections on $\Sigma - p$,
whose holonomy around $p$ has trace $n \exp(2 \pi \sqrt{-1} \frac{d}{n})$. Hence if $\phi$ acts
by the identity on $M$, it will also act by the identity on an
open neighbourhood of $M$ in $\M$, since it acts holomorphically on $\M$. But
since $\M$ is connected, $\phi$ must act by the identity on the entire
$SL(n,\bC)$-moduli space $\M$. Now the generalized \Teim space $\tilde{\T}_p$
of $\Sigma - p$ is also included in $\M$,
hence we get that $\phi$ acts by the identity on $\tilde{\T}_p$. But then the statement
about $\phi$ follows by classical theory of how $\Gamma$ acts on
$\tilde{\T}_p$.
\eproof

\end{document}